\newcount\secno
\newcount\prmno
\newif\ifnotfound
\newif\iffound

% Les commandes
\def\section#1{\vskip1truecm
               \global\def\currenvir{section}
               \global\advance\secno by1\global\prmno=0
               {\bf \number\secno. {#1}}
               \vglue3pt plus 1pt minus 1pt}

\def\subsection{\global\def\currenvir{subsection}
                \global\advance\prmno by1
                \ind{ (\number\secno.\number\prmno) }}
\def\subsec{\global\def\currenvir{subsec}
                \global\advance\prmno by1
                {(\number\secno.\number\prmno)\ }}

\def\proclaim#1{\global\advance\prmno by 1
                {\bf #1 \the\secno.\the\prmno$.-$ }}

\long\def\th#1 \enonce#2\endth{%
   \medbreak\proclaim{#1}{\it #2}\global\def\currenvir{th}\smallskip}

\def\rem#1{\global\advance\prmno by 1
{\it #1} \the\secno.\the\prmno$.-$ }

\magnification 1250
\pretolerance=500 \tolerance=1000  \brokenpenalty=5000
\mathcode`A="7041 \mathcode`B="7042 \mathcode`C="7043
\mathcode`D="7044 \mathcode`E="7045 \mathcode`F="7046
\mathcode`G="7047 \mathcode`H="7048 \mathcode`I="7049
\mathcode`J="704A \mathcode`K="704B \mathcode`L="704C
\mathcode`M="704D \mathcode`N="704E \mathcode`O="704F
\mathcode`P="7050 \mathcode`Q="7051 \mathcode`R="7052
\mathcode`S="7053 \mathcode`T="7054 \mathcode`U="7055
\mathcode`V="7056 \mathcode`W="7057 \mathcode`X="7058
\mathcode`Y="7059 \mathcode`Z="705A
\def\spacedmath#1{\def\packedmath##1${\bgroup\mathsurround =0pt##1\egroup$}
\mathsurround#1
\everymath={\packedmath}\everydisplay={\mathsurround=0pt}}
\def\nospacedmath{\mathsurround=0pt
\everymath={}\everydisplay={} } \spacedmath{2pt}
\def\qfl#1{\buildrel {#1}\over {\longrightarrow}}
\def\phfl#1#2{\normalbaselines{\baselineskip=0pt
\lineskip=10truept\lineskiplimit=1truept}\nospacedmath\smash {\mathop{\hbox to
8truemm{\rightarrowfill}}
\limits^{\scriptstyle#1}_{\scriptstyle#2}}}
\def\hfl#1#2{\normalbaselines{\baselineskip=0truept
\lineskip=10truept\lineskiplimit=1truept}\nospacedmath\smash{\mathop{\hbox to
12truemm{\rightarrowfill}}\limits^{\scriptstyle#1}_{\scriptstyle#2}}}
\def\diagram#1{\def\normalbaselines{\baselineskip=0truept
\lineskip=8truept\lineskiplimit=1truept}   \matrix{#1}}
\def\vfl#1#2{\llap{$\scriptstyle#1$}\left\downarrow\vbox to
6truemm{}\right.\rlap{$\scriptstyle#2$}}
\def\sdir_#1{\mathrel{\mathop{\kern0pt\oplus}\limits_{#1}}}
\def\ssdir_#1^#2{\mathrel{\mathop{\kern0pt\oplus}\limits_{#1}^{#2}}}
\def\pprod_#1{\raise
2pt \hbox{$\mathrel{\scriptstyle\mathop{\kern0pt\prod}\limits_{#1}}$}}
\catcode`\@=11
\font\eightrm=cmr8         \font\eighti=cmmi8
\font\eightsy=cmsy8        \font\eightbf=cmbx8
\font\eighttt=cmtt8        \font\eightit=cmti8
\font\eightsl=cmsl8        \font\sixrm=cmr6
\font\sixi=cmmi6           \font\sixsy=cmsy6
\font\sixbf=cmbx6                
\font\eightgoth=eufm10 at 8pt    
\font\sixgoth=eufm7 at 6pt        \font\fivegoth=eufm5
\newfam\gothfam
\def\eightpoint{%
  \textfont0=\eightrm \scriptfont0=\sixrm \scriptscriptfont0=\fiverm
  \def\rm{\fam\z@\eightrm}%
  \textfont1=\eighti  \scriptfont1=\sixi  \scriptscriptfont1=\fivei
  \textfont2=\eightsy \scriptfont2=\sixsy \scriptscriptfont2=\fivesy
  \textfont\gothfam=\eightgoth \scriptfont\gothfam=\sixgoth
  \scriptscriptfont\gothfam=\fivegoth
  \def\goth{\fam\gothfam\eightgoth}%
   \textfont\itfam=\eightit
  \def\it{\fam\itfam\eightit}%
  \textfont\slfam=\eightsl
  \def\sl{\fam\slfam\eightsl}%
  \textfont\bffam=\eightbf \scriptfont\bffam=\sixbf
  \scriptscriptfont\bffam=\fivebf
  \def\bf{\fam\bffam\eightbf}%
  \textfont\ttfam=\eighttt
  \def\tt{\fam\ttfam\eighttt}%
  \abovedisplayskip=9pt plus 3pt minus 9pt
  \belowdisplayskip=\abovedisplayskip
  \abovedisplayshortskip=0pt plus 3pt
  \belowdisplayshortskip=3pt plus 3pt 
  \smallskipamount=2pt plus 1pt minus 1pt
  \medskipamount=4pt plus 2pt minus 1pt
  \bigskipamount=9pt plus 3pt minus 3pt
  \normalbaselineskip=9pt
  \setbox\strutbox=\hbox{\vrule height7pt depth2pt width0pt}%
  \let\bigf@nt=\eightrm     \let\smallf@nt=\sixrm
  \normalbaselines\rm}
\catcode`\@=12
\def\note#1#2{\footnote{\parindent
0.4cm$^#1$}{\vtop{\eightpoint\baselineskip11pt\hsize15.5truecm\noindent #2}}
\parindent 0cm}
\def\pc#1{\tenrm#1\sevenrm}
\def\up#1{\raise 1ex\hbox{\smallf@nt#1}}
\def\tx{\kern-1.5pt -}
\def\cqfd{\kern 3truemm\unskip\penalty 500\vrule height 4pt depth 0pt width
4pt\smallskip }
\def\ind{\par\hskip 0.8truecm\relax}
\def\inde{\par\nobreak\hskip 0.8truecm\relax}
\def\indp{\par\hskip 0.5truecm\relax}
\def\moins{\mathrel{\hbox{\vrule height 3pt depth -2pt width 6pt}}}
\def\rond{\kern 1pt{\scriptstyle\circ}\kern 1pt}
\def\Hom{\mathop{\rm Hom}\nolimits}
\def\Aut{\mathop{\rm Aut}\nolimits}

\def\at{\mathop{\rm at}\nolimits}
\def\rk{\mathop{\rm rk}\nolimits}
\def\iso{\mathrel{\mathop{\kern 0pt\longrightarrow }\limits^{\sim}}}
\font\spe=cmr8
\def\pro{\hbox{\spe \char 5}}

\vsize = 25.4truecm
\hsize = 16truecm
%\hoffset = -.15truecm
\voffset = -.5truecm
\parindent=0cm
\baselineskip15pt
\overfullrule=0pt\null\vskip0.5cm
\centerline{\bf  Complex manifolds with split tangent bundle}
\smallskip
\smallskip \centerline{Arnaud {\pc BEAUVILLE\note{1}{Partially supported by the
European HCM project ``Algebraic Geometry in Europe" (AGE).\vglue3pt}}} 
\vskip1.2cm

{\bf Introduction}
\smallskip
\ind The theme of this note is to investigate when the tangent bundle of a
compact complex manifold $X$ splits as a direct sum of sub-bundles. This
occurs typically  when  the universal covering space $\widetilde{X}$ of $X$
splits as a product $\pprod_{i\in I} U_i$ of manifolds on which the group $\pi
_1(X)$
 acts diagonally (that is, $\pi _1(X)$ acts on each $U_i$ and 
its action on $\widetilde{X}=\pprod_{} U_i$ is the diagonal action
$g.(u_i)=(gu_i)$): the vector bundles\note{2}{Throughout the paper we will
abuse notation and write $T_{U_i}$ instead of $pr_i^*T_{U_i}$.}
$T_{U_i}$  on
$\widetilde{X}$ are stable under $\pi _1(X)$, hence the decomposition
$T_{\widetilde{X}}=\sdir_iT_{U_i}$ descends to a direct sum decomposition of
$T_X$. For {\it K\"ahler} manifolds, it is tempting to conjecture that the 
converse is true, namely that  any direct sum
decomposition of the tangent bundle $T_X$ (perhaps with the additional 
assumption that the direct summands are integrable) gives rise to a splitting of
the universal  covering. We will show that this is indeed the
case  in three different situations:
\indp {\it a}) $X$ admits a K\"ahler-Einstein metric;
\indp {\it b}) $T_X$ is a direct sum of line bundles of negative degree;
\indp {\it c}) $X$ is a K\"ahler surface.
\ind Case {\it a}) is a direct consequence of the fact that on a compact
K\"ahler-Einstein manifold, any endomorphism of the tangent bundle is
parallel (this idea appears for instance in [Y], and in a more 
implicit form in [K]).  Case {\it b}) is a slight improvement of a
uniformization result of Simpson [S]. To treat case {\it c}) we use the
classification of surfaces and some simple remarks about connections. The
result in this case is actually an easy consequence of the paper [K-O],
where  the authors  classify surfaces with a holomorphic conformal structure
-- this turns out to be closely related to the question we are studying here.
However we found
 simpler and more enlightening to give an independent proof rather than
extracting  from [K-O] the pieces of information that we need.
\ind  In \S 2 we give a few examples which show that for non-K\"ahler 
manifolds a splitting of the tangent bundle does not necessarily imply a
splitting of the universal covering. 
\section{K\"ahler-Einstein manifolds}
{\bf Theorem A}$.-$
{\it  Let $X$ be a compact complex manifold admitting a K\"ahler-Einstein
metric. Assume that the tangent bundle of $X$ has a decomposition
$T_X=\sdir_{i\in I}E_i$. Then the universal covering space of $X$ is a
product $\pprod_{i\in I}U_i$ of complex manifolds, in such a way that the 
decomposition $T_X=\sdir_{i\in I}E_i$ lifts to the decomposition
$T_{\pro U_i}=\sdir_{i\in I}T_{U_i}$; the group $\pi _1(X)$ acts diagonally on
$\pprod_{i\in I} U_i$}.
\smallskip 
	\ind The  proof follows closely that of thm.\ 2.1 in [Y]
(I~am indebted to J. Wahl for pointing out this
reference).\par
{\it Proof}:  (1.1) As a consequence of the Bochner formula,
every endomorphism of $T_X$ is parallel [K]. This applies in
particular to the projectors associated to the direct sum
decomposition of $T_M$; therefore the sub-bundles $E_i$ are
preserved by the hermitian connection, hence the holonomy
representation of $X$ is the direct sum of a family of
representations corresponding to the $E_i$'s.
 By the De Rham theorem, 
 the universal covering space of $X$ splits as a product  $\pprod_{i\in
I}U_i$, such that the  decomposition $T_X=\sdir_{i\in I}E_i$ pulls back
to the decomposition  $T_{\pro U_i}=\sdir_{i\in I}T_{U_i}$.
\ind (1.2) The last assertion  follows from the following
simple observation: {\it if a group $\Gamma $ acting on a product 
$\pprod_{i\in I}U_i$ preserves the  decomposition $T_{\pro U_i}=\sdir_{i\in
I}T_{U_i}
$, it acts diagonally.} Let indeed $\gamma$ be an automorphism of 
$\pprod_{}U_i$; for $j\in I$, put $\gamma_j=pr_j\rond \gamma$. 
The condition $\gamma ^*T_{U_j}=T_{U_j}$  means that the
partial derivatives of $\gamma_j$ in the directions of $U_k$ for $k\not=j$
vanish, hence $\gamma_j((u_i)_{i\in I})$ depends only on $u_j$, which gives 
our claim.\cqfd

\section {Non-K\"ahler examples}
\ind In this section we give examples of manifolds for which the tangent bundle
is a direct sum of line bundles, but which do not satisfy the conclusions of
Theorem A. 
\smallskip 
\subsec{\it Hopf manifolds}
\ind Let $T={\rm diag}(\alpha_1,\dots,\alpha_n)$ be a diagonal matrix, with
$0<|\alpha_i|<1$ for each $i$.  The cyclic group $T^{\bf Z}$ generated by $T$ 
acts freely and properly on ${\bf C}^n\moins\{0\}$; the quotient $X$ is a
compact complex manifold, called a Hopf manifold.
For each
non-zero complex number $\theta$,  denote by $L_\theta$ the
flat line bundle associated to the  character of $\pi _1(X)=T^{\bf Z}$
mapping $T$ to $\theta$; in other words, $L_\theta$ is the
quotient of the trivial line bundle   $({\bf C}^n\moins\{0\})\times {\bf C}$
by the action of  the automorphism $(T,\theta)$. By construction we have
$T_X=\ssdir_{i=1}^{n}L_{\alpha_i}$, but  the universal covering space
${\bf C}^n\moins\{0\}$ of $X$ is clearly not a product.
\smallskip
\subsec{\it Complex compact nilmanifolds}\par\nobreak
\hskip0.8truecm These are compact manifolds  $X=G/\Gamma $, where $G$ is a
 nilpotent complex Lie group and $\Gamma $ a
discrete subgroup of $G$. We may assume that $G$ is simply-connected and
non-commutative (to exclude the trivial case of complex tori). A well-known
example is the Iwasawa manifold $U({\bf C})/U({\bf Z}[i])$, where
$U$ is  the group of upper-triangular $3\times 3$ matrices with diagonal
entries $1$; many examples can be obtained in an analogous way. 
\ind The tangent bundle of  $X=G/\Gamma $ is trivial, and its
universal covering space $G$  is isomorphic to ${\bf C}^n$; however we claim
that whatever isomorphism $G\iso {\bf C}^n$ we choose, {\it the action
of $\Gamma $ cannot be diagonal}. Indeed if $\Gamma $ acts diagonally, the
standard trivialization of $T_{{\bf C}^n}$ deduced from
the coordinate system descends to a  trivialization   of
$T_X$. Any such trivialization lifts to a trivialization of $T_G$   
defined by a basis of right invariant vector fields; therefore the standard
trivialization of $T_{{\bf C}^n}$ is $G$\tx equivariant. In view of 1.2 this
means that $G$ itself acts diagonally on ${\bf C}^n$, hence $G$ embeds into
$\Aut({\bf C})^n$. Now any nilpotent connected subgroup of the affine group
$\Aut({\bf C})$ is commutative, so we conclude that $G$ is 
 commutative, contrary to our hypothesis.

\section{Simpson's uniformization result}
\ind The following lemma, which is a variation on the Baum-Bott theorem
[B-B], will allow us to slightly improve  Simpson's  result:
\th Lemma
\enonce Let $X$ be a complex manifold, and $E$ a direct summand of $T_X$. The
Atiyah class $\at(E)\in H^1(X,\Omega^1_X\otimes {\cal E}nd(E))$ 
 comes from $H^1(X,E^*\otimes {\cal E}nd(E))$. In particular, any class in
$H^r(X,\Omega_X^r)$ given by a polynomial in the Chern classes of $E$
vanishes for $r>\rk(E)$.
\endth
{\it Proof}: Write $T_X=E\oplus F$; let $p:T_X\rightarrow E$ be the
corresponding projection. For any 
sections $U$ of $E$ and  $V$ of $F$ over some open subset of $X$, put
$D_VU=$ $p([V,U])$. This expression is ${\cal O}_X$\tx linear in $V$ and
satisfies the Leibnitz rule $D_V(fU)=fD_V(U)+(Vf)U$,  so that $D$ is a
$F$\tx {\it connection}  on $E$ [B-B]: 
 if we denote by ${\cal D}^1(E)$  the sheaf of
differential operators  $\Delta:E\rightarrow E$, of degree $\leq 1$, whose
symbol $\sigma(\Delta)$ is scalar, this means that $D$ defines an ${\cal
O}_X$\tx linear map
$F\rightarrow {\cal D}^1(E)$ such that $\sigma(D_V)=V$ for all local sections
$V$ of $F$. Thus the exact sequence
$$ 0 \rightarrow   {\cal E}nd(E)  \longrightarrow   {\cal D}^1(E)
\qfl{\sigma}  T_X\rightarrow 0
$$
splits over the sub-bundle $F\i T_X$; therefore its extension class
$\at(E)\in$\break $H^1(X,\Omega^1_X\otimes {\cal E}nd(E))$ vanishes
in $H^1(X,F^*\otimes {\cal E}nd(E))$, hence comes from\break $H^1(X,E^*\otimes
{\cal E}nd(E))$. The last assertion follows from the definition of the Chern
classes in terms of the Atiyah class.\cqfd\smallskip
\ind We denote as usual by ${\bf H}$ the Poincar\'e upper half-space.\par

{\bf Theorem B}$.-$ {\it Let $X$ be a compact  K\"ahler
manifold, with a K\"ahler class $\omega $. Assume that the tangent bundle
$T_X$ is a direct sum of  line bundles
$L_1,\ldots,L_n$ with $\omega ^{n-1}.\, c_1(L_i)<0$ for each $i$.
Then the
 universal covering space of $X$ is ${\bf H}^n$, and the 
decomposition  $T_X=\oplus L_i$ lifts to the canonical decomposition}
$T_{{\bf H}^n}=(T_{\bf H})^{\oplus n}$.
\smallskip
 {\it Proof}: This is Cor.\
9.7 of [S], except that Simpson makes the extra hypothesis $\omega
^{n-2}.\,(c_1(X)^2-2c_2(X))=0$ (the assertion about the compatibility of
decompositions is not stated in {\it loc. cit.}, but follows directly from
the proof). Now lemma 3.1 gives $c_1(L_i)^2=0$ for each $i$, hence
$c_1(X)^2-2c_2(X)=0$.\cqfd

\section{The surface case}
{\bf Theorem C}$.-$
{\it   Let $X$ be a compact complex surface. The tangent
bundle of $X$ splits as a direct sum of two line bundles if and only if one of
the following occurs:
\indp {\rm a)} The universal covering space of $X$ is a product $U\times V$ of
two (simply-connected) Riemann surfaces and the group $\pi_1(X)$ acts
diagonally on $U\times V$; in that case the given splitting of $T_X$ lifts to
the direct sum decomposition  $T_{U\times V}=T_U\oplus T_V$.
\indp {\rm b)} $X$ is a Hopf surface, with
universal covering space  ${\bf
C}^2\moins\{0\}$. Its fundamental group is isomorphic to ${\bf
Z}\oplus{\bf Z}/m{\bf Z}$, for some integer $m\ge 1$; it is generated by a
diagonal automorphism
$(x,y)\mapsto (\alpha x,\beta y)$  with $|\alpha |\le|\beta |<1$, and a 
 diagonal automorphism $(x,y)\mapsto (\lambda x,\mu y)$ where $\lambda $
and $\mu $ are primitive $m$\tx th roots of $1$}.
\smallskip  \ind  As a corollary,
for K\"ahler surfaces we see  that any direct sum decomposition of the 
tangent bundle
gives rise to a splitting of the universal  covering, as announced in the
introduction. \smallskip 
\subsection Before starting the proof we will need a few preliminaries.
From now on we denote by  $X$ a compact
complex surface;
 we assume given a direct sum decomposition $\Omega^1_X\cong L\oplus M$.  By
lemma 3.1 (or by [B-B])  the Chern  class
$c_1(L) \in H^1(X,\Omega ^1_X)$ belongs to the subspace $H^1(X,L)$, and
similarly for $M$.  As a consequence, we get:
\subsection We have $L^2=M^2=0$, and therefore
$c_1^2(X)=2\,L.M=2c_2(X)$. 
\smallskip 
\ind The following consequence is less obvious.
\th Proposition
\enonce Let $C$ be a smooth rational curve in $X$. Then $C^2\ge 0$.
\endth
{\it Proof}: Put $C^2=-d$ and assume $d>0$. Since $H^1(C,{\cal O}_C(d+2))=0$,  
the exact sequence
$$0\rightarrow {\cal O}_C(d)\longrightarrow \Omega ^1_{X\,|C}\longrightarrow
\Omega ^1_C\rightarrow 0$$splits, giving an isomorphism $\Omega
^1_{X\,|C}\cong  {\cal O}_C(d)\oplus {\cal O}_C(-2)$. Thus one of the line
bundles $L$ or $M$, say $L$, satisfies $L_{|C}\cong {\cal O}_C(d)$. Consider
the commutative diagram
\vskip-10pt$$\diagram{
H^1(X,L) & \hfl{}{} & H^1(X,\Omega ^1_X)&\cr
\vfl{}{} & &\vfl{}{} &\cr
H^1(C,L_{|C}) & \hfl{}{} & H^1(C,\Omega ^1_C)&;\cr
}$$\vskip-10pt
since $d>0$ we have $H^1(C,L_{|C})=0$; thus $c_1(L)$ goes to $0$ in
$H^1(C,\Omega ^1_C)$,  which means $d=0$, a contradiction.\cqfd

\subsection We shall come across  situations where the vector bundle
$\Omega ^1_X=L\oplus M$ appears as an extension
$$0\rightarrow P\longrightarrow  \Omega ^1_X\qfl{p} Q\rightarrow 0$$
of two line bundles $P$ and $Q$. In that case, 
	\indp -- either the restriction of $p$ to one of the direct summands of
$\Omega^1_X$, say
$M$, is surjective; then the exact sequence splits, $Q$ is isomorphic to $M$ 
and $P$ to $L$;
\indp -- or the restriction of $p$ to both $L$ and  $M$ is not surjective; then
there exists effective (non-zero) divisors $A$ and $B$, whose supports do not
intersect, such that 
$L\cong Q(-A)$, $M\cong Q(-B)$ and $P\cong Q(-A-B)$; the exact
sequence does {\it not} split. \ind In particular, if
$\Hom(P,Q)=0$, the exact sequence splits.
\smallskip
\subsection Finally we will need  some classical facts about connections 
(see [E]). Let
$p:M\rightarrow B$ be a smooth holomorphic map between complex manifolds, whose
fibres are isomorphic to a fixed variety $F$.  A
{\it connection} on
$p$ is a splitting of the exact sequence
$$0\rightarrow p^*\Omega^1_B\longrightarrow \Omega^1_M\longrightarrow
\Omega^1_{M/B}\rightarrow 0\ ,$$that is a sub-bundle $L\i \Omega^1_M$ mapping
isomorphically onto
$\Omega^1_{M/B}$; the connection is flat (or integrable) if $dL\i L\wedge
\Omega^1_M$  (this  is automatic if $B$ is a curve). In that case  the
group $\pi _1(B)$ acts on $F$ by complex automorphisms, and $M$ is the fibre 
bundle on $B$ with fibre $F$  associated to the universal covering 
$\widetilde{B}\rightarrow B$, that is the quotient of
$\widetilde{B}\times F$ by the group  $\pi _1(B)$ acting diagonally; the
splitting $\Omega ^1_M=p^*\Omega ^1_B\oplus L$ pulls back to the
decomposition $\Omega ^1_{\widetilde{B}\times F}=\Omega
^1_{\widetilde{B}}\oplus \Omega^1_F$.

\section{Proof of theorem C}
\subsec{\it Kodaira dimension $2$}
\inde If $\kappa (X)=2$, the canonical bundle $K_X$ 
is ample by Prop.\ 4.3. The Aubin-Calabi-Yau theorem
 implies that $X$ admits a K\"ahler-Einstein metric; we can therefore apply 
Theorem A.

\medskip 
\subsec{\it Kodaira dimension $1$} 
\inde If $\kappa (X)=1$, $X$ admits an
elliptic fibration $p:X\rightarrow B$. By 4.2  we have $c_2(X)=0$;
this implies that the only singular fibres of $p$ are multiples of smooth
elliptic curves (see [B1], VI.4 and VI.5). For $b\in B$, we write
$p^*[b]=m(b)\,F_b$, where $F_b$ is a smooth elliptic curve; we have $m(b)\ge
1$ and $m(b)=1$ except for finitely many points. Put $\Delta =\sum_b
(m(b)-1)\,F_b$. We have an exact sequence $$0\rightarrow p^*\Omega^1_B(\Delta
)\longrightarrow \Omega ^1_X\longrightarrow \omega_{X/B}\rightarrow 0\
,\leqno\subsec$$ where $\omega_{X/B}$ is the relative dualizing line bundle.
Since $\chi({\cal O}_X)=0$ by Riemann-Roch, we deduce from [B-P-V], V.12.2 and
III.18.2, that  $\omega_{X/B}$ is a torsion line bundle. Since
$K_X=p^*\Omega^1_B(\Delta )\otimes \omega_{X/B}$, the hypothesis $\kappa
(X)=1$ implies $\Hom(p^*\Omega^1_B(\Delta), \omega_{X/B})=0$, hence 
 the exact sequence (5.3) splits by  4.4.  \ind Let $\rho
:\widetilde{B}\rightarrow B$ be the orbifold universal covering of $(B, m)$:
this is a ramified Galois covering,  with $\widetilde{B}$ simply-connected, 
such  that the stabilizer of a point $\tilde b\in \widetilde{B}$ is a cyclic
group of order $m(\rho (\tilde b))$ (see for instance [K-O], lemma 6.1; note 
that because of the hypothesis  $\kappa (X)=1$ and the formula for $K_X$,
there are at least $3$ multiple fibers if $B$ is of genus $0$). Let   
$\widetilde{X}$ be the normalization of $X\times
_B\widetilde{B}$. We have a commutative diagram
$$\diagram{\widetilde{X} & \hfl{\pi }{} & X\cr
\vfl{\tilde p}{} && \vfl{}{p} \cr
\widetilde{B} & \hfl{\rho}{} & B\cr
}$$where $\tilde p$ is smooth and $\pi $ is \'etale $([B1],\, VI.7')$. The 
exact sequence $$0\rightarrow \tilde
p^*\Omega^1_{\widetilde{B}}\longrightarrow
\Omega^1_{\widetilde{X}}\longrightarrow
\Omega^1_{\widetilde{X}/\widetilde{B}}\rightarrow 0 $$ coincides with the
pull back under $\pi$ of the exact sequence  (5.3);
therefore $p$ admits an integrable connection, given by the subbundle  $\pi
^*M $ of $\Omega^1_{\widetilde{X}}$. The result follows from 4.5 and
1.2.
\medskip
\subsec{\it Kodaira dimension $0$}
\inde Assume $\kappa (X)=0$. By 4.2  and the classification of
surfaces, $X$ is either a complex torus, a biellip\-tic surface, or a
Kodaira surface. Complex tori and biellip\-tic surfaces fall into case {\it
a}) of the theorem (a bielliptic surface is the quotient of a product
$E\times F$ of elliptic curves by a finite abelian group acting diagonally).
\ind A primary Kodaira surface has trivial canonical bundle and admits a
smooth elliptic fibration $p:X\rightarrow B$. Thus the exact sequence (4.2)
realizes $\Omega ^1_X$ as an extension of ${\cal O}_X$ by ${\cal O}_X$. Since
$h^{1,0}(X)=1$, this extension is non-trivial, and it follows from 4.4
that $\Omega ^1_X$ does not split.
\ind A secondary Kodaira surface admits a primary Kodaira surface as a finite
\'etale cover, hence its tangent bundle cannot split either.
\medskip
\subsec{\it Ruled surfaces}
\ind  We  consider the case when $X$ is algebraic and  $\kappa (X)=-\infty$.
By 4.2  and 4.3, $X$ is a geometrically ruled surface, that is a
projective bundle $p:X\rightarrow B$ over a curve.  We  again consider the
exact sequence $$0\rightarrow p^*\Omega^1_B\longrightarrow
\Omega ^1_X\longrightarrow \Omega^1_{X/B}\rightarrow 0\ ;$$
since $\Omega^1_{X/B}$ has negative degree on the fibres, we have
$\Hom(p^*\Omega^1_B,\Omega^1_{X/B})=0$, hence  by
4.4 the above exact sequence splits: one of
the direct summands of $\Omega ^1_X$ defines an integrable connection for
$p$. The result follows then from 4.5.
\medskip
\subsec{\it Inoue surfaces}
\ind  We now assume
 that $X$  is not algebraic and $\kappa(X) =-\infty$, so that $X$ is what
is usually called a surface  of type$\,VII_0$. These surfaces have
$b_1=h^{0,1}=1$
 and therefore $c_1^2+c_2=12\chi ({\cal O}_X)=0$; in
our case this gives $c_1^2=c_2=0$ in view of 4.2, and finally $b_2=0$.
Moreover we have $H^0(X,\Omega^1_X\otimes L^{-1}) \not=0$. The surfaces with
these properties have been completely classified by Inoue [I]: they are either
  Hopf surfaces, or belong to three classes of surfaces constructed by Inoue
({\it loc.\ cit.}).
\ind We first consider the Inoue surfaces. The surfaces
$S_M$ of the first class are quotients of ${\bf H}\times {\bf C}$ by a group
acting diagonally, hence they fall into case {\it a}) of the
theorem.  \ind The surfaces $S_{N,p,q,r;t}^{(+)}$ of the
second class are quotients of ${\bf H}\times {\bf
C}$ by a group which does {\it not} act diagonally. This
action leaves invariant the vector field $\partial/\partial
z$ on ${\bf C}$, which therefore descends to a non-vanishing
vector field $v$ on $X$. This gives rise to an exact sequence
$$0\rightarrow K_X\ \qfl{i(v)}\ \Omega ^1_X\ \qfl{i(v)}\
{\cal O}_X\rightarrow 0\ , $$  which does not split since
$h^{1,0}(X)=0$. We have $H^0(X,K_X^{-1})=0$, for instance
because $X$ contains no curves; we infer from 4.4 that
$\Omega ^1_X$ does not split. \ind The surfaces $S_{N,p,q,r}^{(-)}$
of the third class are quotients of certain
surfaces of the second class by a fixed point free
involution; therefore their tangent bundle does not split
either.
 
\medskip
\subsec{\it Primary Hopf surfaces} 
\ind   It remains to consider the class of Hopf
surfaces, which are by  definition the
surfaces of class$\,VII_0$ whose universal covering space is ${\bf
W}:={\bf C}^2\moins\{0\}$. We consider first the {\it
primary} Hopf surfaces, which  are quotients of ${\bf W}$ by
the infinite cyclic group generated by an automorphism $T$ of
${\bf W}$. According to [Ko], \S 10, there are two cases to
consider: \indp {\it a}) $T(x,y)=(\alpha x,\beta y)$ for some
complex numbers $\alpha,\beta$ with $0<|\alpha|\le
|\beta|<1$; \indp {\it b}) $T(x,y)=(\alpha^m x+\lambda
y^m,\alpha y)$ for some positive integer $m$ and non-zero complex
numbers $\alpha,\lambda$ with $|\alpha|<1$. \ind As in 2.1, we denote by
$L_\theta$, for $\theta \in{\bf C}$, the flat line bundle associated to the 
character of $\pi _1(X)$ mapping $T$ to $\theta$.
In case {\it a}) we find $\Omega^1_X=L_\alpha^{-1} \oplus L_\beta^{-1} $,
so the tangent bundle splits. \ind Let us consider  case {\it
b}).  The form $dy$ on ${\bf W}$ satisfies $T^*dy=\alpha\,dy$, hence  
descends to a form $\overline{dy}$ in $H^0(X,\Omega_X^1\otimes L_\alpha)$;
similarly the function $y$  descends to a non-zero section  of $L_\alpha $. We
have an  exact sequence
$$0\rightarrow L_\alpha ^{-1}\ \qfl{\overline{dy}}\ \Omega ^1_X\
\longrightarrow \ L_\alpha ^{-m}\rightarrow 0\ .$$ 
\ind Since $L_\alpha$ has a non-zero section, the space
$\Hom(L_\alpha ^{-1},L_\alpha ^{-m})$ is zero for $m>1$. Hence if $\Omega_X^1$
splits, we deduce from 4.4 that the exact sequence splits.
 This means that there exists  a form $\overline{\omega} \in H^0(X,\Omega_X^1
\otimes L_\alpha ^m)$ such that $\overline{\omega} \wedge \overline{dy}\not=0$.
Then $\overline{\omega} \wedge \overline{dy}$ is a generator of the trivial 
line bundle $K_X\otimes L_\alpha^{m+1}$, hence pulls back to $c\,dx\wedge dy$ 
on  ${\bf W}$, for some constant $c\not=0$. Therefore the pull back $\omega$ 
of $\overline{\omega}$ to ${\bf W}$ is of the form  $c\,dx+f(x,y)dy$ for some
holomorphic function $f$ on ${\bf C}^2$. The flat line bundle  $L_\alpha^m$
carries a flat holomorphic connection  $\nabla$;  the 2-form $\nabla
\overline{\omega}$, which is a global section of  $K_X\otimes L_\alpha^m\cong
L_\alpha^{-1}$,  is zero. This implies  $d\omega =0$, so the function $f(x,y)$
is independent of $x$; let us write it  $f(y)$. Now the condition $T^*\omega
=\alpha^m\omega $ reads $\alpha f(\alpha y)+c\lambda  my^{m-1}=\alpha ^mf(y)$.
Differentiating $m$ times we find  $f^{(m)}=0$, then differentiating $m-1$
times leads to a contradiction.
\medskip
\subsec{\it Secondary Hopf surfaces}
\ind A secondary Hopf surface $X$ is the quotient of ${\bf W}$ by a
group $\Gamma $ acting freely, containing  a central, finite
index subgroup generated by an automorphism $T$ of the above type. We assume
that $\Omega ^1_X$ splits. The primary Hopf surface $Y={\bf W}/
T^{\bf Z}$ is a finite \'etale cover of $X$, so $\Omega ^1_Y$ also splits; it
follows from 5.7 that  $T$ is of type {\it a}), and  that  $\Gamma $ does
not contain any transformation of type {\it b}). According to [Ka], \S 3, this
implies that after an appropriate change
of coordinates, the group $\Gamma $ acts {\it linearly} on ${\bf C}^2$. 
\ind We claim that $\Gamma $ is contained in a maximal torus of ${\bf
GL}(2,{\bf C})$. This is clear if $\alpha \not=\beta $, because  $T$ is
central in $\Gamma $.
If $\alpha =\beta $, the direct sum decomposition of  $\Omega ^1_X$ pulls back
to a decomposition $\Omega ^1_Y=L_\alpha^{-1} \oplus L_\alpha^{-1} $ (5.7),
which for an appropriate choice of coordinates comes from the decomposition
$\Omega ^1_{\bf W}={\cal O}_{\bf W}^{}dx\oplus {\cal O}_{\bf W}^{}dy$. Since
$\Gamma $ must preserve this decomposition, it is  contained in the diagonal
torus. \ind Thus we may identify $\Gamma $ with a subgroup of $({\bf C}^*)^2$;
since it acts freely on ${\bf W}$, the first projection $\Gamma \rightarrow
{\bf C}^*$ is injective. Therefore the torsion subgroup  of $\Gamma $ is
cyclic, and we are in case {\it b}) of the theorem.\cqfd

\vskip1.5cm
\centerline{ REFERENCES} \vglue15pt\baselineskip12pt
\def\num#1{\smallskip\item{\hbox to\parindent{\enskip [#1]\hfill}}}
\parindent=1.3cm 
\num{B1} A.\ {\pc BEAUVILLE}: {\sl 	Surfaces alg\'ebriques complexes}.  
Ast\'e\-risque {\bf 54} (1978).
\num{B2} A.\ {\pc BEAUVILLE}: {\sl Vari\'et\'es k\"ahl\'eriennes dont la 
premi\`ere classe de Chern est nulle}.  J.\ of Diff.\ Geo\-metry {\bf 18},
755--782  (1983).  
\num{B-B} P.\ {\pc BAUM}, R.\ {\pc BOTT}: {\sl On the zeros
of meromorphic vector-fields}.  Essays on Topology and Related Topics, 
29--47; Springer, New York (1970). 
\num{B-P-V} W.\ {\pc BARTH}, C.\ {\pc PETERS}, A.\ {\pc VAN DE} {\pc VEN}:
{\sl Compact complex surfaces}. Ergebnisse der Math., Springer-Verlag (1984).
\num{E} C.\ {\pc EHRESMANN}: {\sl Les connexions infinit\'esimales dans un 
espace fibr\'e diff\'eren\-tiable}.  Colloque de topologie, Bruxelles (1950),
 29--55. G.\ Thone, Li\`ege (1951).
\num{I} M.\ {\pc INOUE}: {\sl On surfaces of class$\,VII_0$}. Invent.\ math.\ 
{\bf 24}, 269--310 (1974). 
\num{Ka} M.\ {\pc KATO}: {\sl Topology of Hopf surfaces}. J.\ Math.\ Soc.\ 
Japan {\bf 27}, 222--238 (1975). 
\num{K} S.\ {\pc KOBAYASHI}: {\sl
First Chern class and holomorphic tensor
fields}. Nagoya Math. J. {\bf 77}, 5--11 (1980).
\num{K-O} S.\ {\pc KOBAYASHI}, T.\ {\pc OCHIAI}: {\sl Holomorphic structures 
modeled after hyperquadrics}. T\^ohoku Math.\ J.\ {\bf 34},
587--629  (1982).
\num{Ko} K.\ {\pc KODAIRA}: {\sl On the structure of compact complex analytic
surfaces} II. Amer.\  J.\  of Math.\  {\bf 88}, 682--721 (1966).
\num{S} C.\ {\pc SIMPSON}: {\sl Constructing variations of Hodge structure
using Yang-Mills theory and applications to uniformization}. J. Amer. Math. Soc.
{\bf 1}, 867--918  (1988). 
\num{Y} S.-T.\ {\pc YAU}: {\sl A splitting theorem and an algebraic
geometric characterization of locally hermitian symmetric spaces}. Comm. in
Analysis and Geometry {\bf 1}, 473--486  (1993).
\vskip1cm
\def\pc#1{\eightrm#1\sixrm}
\hfill\vtop{\eightrm\hbox to 5cm{\hfill Arnaud {\pc BEAUVILLE}\hfill}
 \hbox to 5cm{\hfill DMI -- \'Ecole Normale
Sup\'erieure\hfill} \hbox to 5cm{\hfill (URA 762 du CNRS)\hfill}
\hbox to 5cm{\hfill  45 rue d'Ulm\hfill}
\hbox to 5cm{\hfill F-75230 {\pc PARIS} Cedex 05\hfill}}
\end